\documentclass[journal]{new-aiaa}

\usepackage[utf8]{inputenc}
\usepackage{textcomp}
\usepackage{graphicx}
\usepackage{amsmath,mathtools}
\usepackage{siunitx}
\usepackage{longtable,tabularx}
\usepackage{booktabs}
\usepackage{xcolor}

\graphicspath{{figures/}}
\hypersetup{hidelinks}

\newcommand{\vect}[1]{\boldsymbol{#1}}
\newcommand{\norm}[1]{\left\lVert #1\right\rVert}
\newcommand{\prox}{\operatorname{prox}}
\newcommand{\diag}{\operatorname{diag}}
\newcommand{\ind}{\mathcal I}
\newcommand{\R}{\mathbb R}
\newcommand{\zero}{\vect 0}

\title{Microsecond-Class Powered-Descent Optimization via Exact Condensation and Strong Convex Regularization}

\author{Wenbo Li\textsuperscript{1}, Ziqi Xu\textsuperscript{2},
Dai Shen\textsuperscript{2}, and Shengping Gong\textsuperscript{2,3,*}}
\affil{\textsuperscript{1}School of Aerospace Engineering, Tsinghua University,
Beijing 100084, China\\
\textsuperscript{2}School of Astronautics, Beihang University,
Beijing 102206, China\\
\textsuperscript{3}State Key Laboratory of High-Efficiency Reusable
Aerospace Transportation Technology, Beijing 102206, China\\
\textsuperscript{*}Corresponding author}

\begin{document}

\maketitle
\begin{center}
\small Manuscript prepared July 27, 2026
\end{center}

\begin{abstract}
Fuel-dominant powered descent can be written as a convex program, but the usual full-state epigraph formulation still carries many state variables, fuel epigraph variables, and dynamics equalities. In addition, the pure-fuel objective provides no strong-convexity curvature. This paper combines three structural reductions. First, a dimensionally consistent low-weight energy term makes the control solution unique. Second, a terminal-state sensitivity recursion eliminates every intermediate state exactly; invertible row normalization turns the 30-node baseline with 300 primal variables and 174 equality multipliers into a problem with 90 control variables and six terminal multipliers. Third, the shared radial structure of the fuel norm and thrust ball gives an exact closed-form proximal operator consisting of group shrinkage followed by magnitude clipping. The condensed problem is solved with a fixed-budget extrapolated proportional--integral projected-gradient iteration implemented in fixed-size C17 arrays. In a Mars powered-descent case with energy weight 0.02 and relative reference-solution tolerance $10^{-3}$, the iteration count decreases from 2744 for the pure-fuel full-state epigraph baseline to 93, while the fuel metric increases by only 0.033\%. The mean end-to-end solve time is \SI{68.2}{\micro\second}, and P99 is \SI{128.1}{\micro\second}, on an Intel i7-10875H. Because the thrust set in this test case is already a convex ball, the contribution is fast solution of the convex core rather than a new lossless-convexification theorem.
\end{abstract}

\section*{Nomenclature}

{\renewcommand\arraystretch{1.0}
\noindent\begin{longtable*}{@{}l @{\quad=\quad} p{0.79\textwidth}@{}}
$A,B_0,B_1$ & discrete state-transition and control matrices \\
$D$ & invertible row-scaling matrix \\
$G,H_c$ & unscaled and row-normalized condensed dynamics matrices \\
$J_E,J_F$ & energy and fuel metrics \\
$N$ & number of control nodes \\
$P$ & block-diagonal quadratic-objective matrix \\
$q_k$ & nodewise fuel proximal coefficient \\
$\vect r,\vect v,\vect u$ & position, velocity, and thrust acceleration \\
$t_f,\Delta t$ & final time and uniform time step \\
$\vect U$ & stacked normalized control vector \\
$u_{\max}$ & maximum thrust acceleration \\
$\alpha,\beta$ & primal and dual stepsizes \\
$\lambda$ & energy weight in the fuel-dominant objective \\
$\mu$ & strong-convexity modulus \\
$\rho$ & extrapolation parameter \\
$\vect\eta$ & terminal-equality multiplier \\
\multicolumn{2}{@{}l}{Superscripts}\\
$(\cdot)^\star$ & reference optimum \\
$(\hat{\cdot})$ & normalized variable \\
\end{longtable*}}
\addtocounter{table}{-1}

\section{Introduction}

\lettrine{P}{owered-descent} guidance must generate a feasible thrust history within a short guidance cycle while satisfying terminal position, terminal velocity, and actuator constraints. Convex optimization provides globally optimal solutions for important translational landing models \cite{acikmese2007convex,blackmore2010minimum}. Lossless convexification extends this capability to nonconvex thrust lower bounds under stated regularity assumptions, and successive convexification treats more general six-degree-of-freedom dynamics and state-triggered constraints \cite{szmuk2020successive,malyuta2021advances}. These developments have made optimization-based guidance credible, but onboard use is governed by more than mathematical solvability. Memory, execution-path regularity, latency dispersion, and the ability to certify the implementation can be as important as nominal solution accuracy.

Interior-point methods solve conic programs to high accuracy in relatively few iterations, but each iteration requires the formation and factorization of a Karush--Kuhn--Tucker system. Code generation, sparse factorization, and embedded conic solvers reduce this burden substantially \cite{mattingley2012cvxgen,domahidi2013ecos,stellato2020osqp}. First-order primal--dual methods offer a complementary route: they replace factorization with matrix--vector products and projections. The proportional--integral projected-gradient method provides convergence guarantees for conic optimization, and its extrapolated form can improve numerical convergence \cite{yu2022pipg,yu2023xpipg}. Customized variants have demonstrated the relevance of this approach to powered-descent guidance \cite{elango2022customized,kamath2023customized}.

The computational potential of a first-order method is not determined by the iteration formula alone. A standard fuel epigraph retains all discrete states, introduces one scalar epigraph variable per node, and produces a long chain of dynamics multipliers. Its objective is convex but not strongly convex, so control directions can drift along nearly flat optimal faces. In contrast, a terminal-only translational landing model does not require intermediate states as online decision variables. Moreover, the Euclidean fuel norm and a concentric thrust ball admit a joint closed-form proximal map. These two structural facts suggest that formulation design should precede low-level code optimization.

This paper develops that formulation-to-code chain. The contribution is not a new lossless-convexification theorem or a new generic primal--dual method. It is the coordinated construction of a small, strongly convex, proximable problem whose complete online iteration can be implemented without factorization or dynamic memory. The specific contributions are:

\begin{enumerate}
\item A dimensionally consistent energy--fuel objective is scaled into a block-separable form with explicit strong-convexity modulus $\mu=\lambda/2$. This curvature guarantees a unique control solution and exposes the fuel-versus-convergence tradeoff.
\item A terminal-sensitivity recursion exactly eliminates every intermediate state and every interval dynamics multiplier. For $N=30$, the resulting problem has 90 primal variables and only six equality multipliers, compared with 300 and 174 in the epigraph baseline.
\item The fuel norm and thrust ball are combined into a single radial proximal subproblem. Its exact solution is group soft thresholding followed by magnitude clipping, with no epigraph variable and no inner iteration.
\item A fixed-budget C17 realization uses static arrays and precomputed spectral stepsizes. Cross-language verification and timing tests show a mean latency of \SI{68.2}{\micro\second} and a 99th-percentile latency of \SI{128.1}{\micro\second} for the fixed test case.
\end{enumerate}

The remainder of the paper derives the condensed problem, proves the strong-convexity and proximal results, presents the fixed-budget iteration, and evaluates the accuracy, tradeoffs, and timing limitations.

\section{Powered-Descent Problem}

\subsection{Continuous Dynamics and Constraints}

Consider a fixed-mass, three-degree-of-freedom translational model. The vertical direction is the $y$ axis. Position, velocity, thrust acceleration, and constant gravitational acceleration are denoted by $\vect r,\vect v,\vect u,\vect g\in\R^3$, respectively:
\begin{equation}
\dot{\vect r}=\vect v,\qquad
\dot{\vect v}=\vect u+\vect g .
\label{eq:dynamics}
\end{equation}
The boundary conditions are
\begin{equation}
\vect r(0)=\vect r_0,\quad \vect v(0)=\vect v_0,\quad
\vect r(t_f)=\vect r_f,\quad \vect v(t_f)=\vect v_f ,
\label{eq:boundary}
\end{equation}
and the actuator envelope is
\begin{equation}
0\leq\norm{\vect u(t)}_2\leq u_{\max}.
\label{eq:thrust}
\end{equation}
The verification model contains no nonzero lower-thrust bound, pointing cone, or state-path constraint. Thus, Eq.~\eqref{eq:thrust} is already a convex Euclidean ball. This distinction matters: the present microsecond result isolates the computational effect of strong convexification, condensation, and proximal structure rather than claiming a new relaxation result for the more general nonconvex powered-descent problem.

\subsection{Trapezoidal Discretization and Pure-Fuel Baseline}

Divide $[0,t_f]$ into $N-1$ uniform intervals with $\Delta t=t_f/(N-1)$. Trapezoidal discretization gives
\begin{align}
\vect v_{k+1}&=\vect v_k+\frac{\Delta t}{2}
(\vect u_k+\vect u_{k+1})+\Delta t\vect g,
\label{eq:trapv}\\
\vect r_{k+1}&=\vect r_k+\frac{\Delta t}{2}
(\vect v_k+\vect v_{k+1}),\qquad k=1,\ldots,N-1 .
\label{eq:trapr}
\end{align}
The discretized fuel metric is
\begin{equation}
J_F=\Delta t\sum_{k=1}^{N}w_k\norm{\vect u_k}_2,\qquad
w_1=w_N=\frac12,\quad w_k=1\ (2\leq k\leq N-1).
\label{eq:fuel}
\end{equation}
A standard second-order-cone epigraph introduces $\sigma_k$ at every node:
\begin{align}
\min_{\vect r_k,\vect v_k,\vect u_k,\sigma_k}\quad&
\Delta t\sum_{k=1}^{N}w_k\sigma_k,
\label{eq:epigraph_obj}\\
\text{s.t.}\quad&
\text{Eqs.~\eqref{eq:boundary}, \eqref{eq:trapv}, and \eqref{eq:trapr}},\nonumber\\
&\norm{\vect u_k}_2\leq\sigma_k\leq u_{\max}.
\label{eq:epigraph_cone}
\end{align}
For $N$ nodes, this form contains $10N$ primal variables. It also retains $6(N-1)$ interval dynamics equalities. The baseline first-order solver operates directly on this full epigraph form.

\subsection{Dimensionless Fuel-Dominant Strongly Convex Objective}

Define an energy metric
\begin{equation}
J_E=\Delta t\sum_{k=1}^{N}w_k\norm{\vect u_k}_2^2
\label{eq:energy}
\end{equation}
and the dimensionless weighted objective
\begin{equation}
J_\lambda=
\lambda\frac{J_E}{t_fu_{\max}^2}
{}+(1-\lambda)\frac{J_F}{t_fu_{\max}},\qquad 0<\lambda\leq1 .
\label{eq:weighted}
\end{equation}
The denominators are the natural scales of the two metrics, so $\lambda$ is not contaminated by physical units. As $\lambda\rightarrow0^+$, Eq.~\eqref{eq:weighted} approaches the pure-fuel objective, whereas every positive $\lambda$ adds curvature. The numerical study uses $\lambda=0.02$ as a fuel-dominant setting and separately reports its fuel penalty relative to the pure-fuel solution.

Introduce diagonal normalizations
\begin{equation}
\vect r=P_r\hat{\vect r},\qquad
\vect v=P_v\hat{\vect v},\qquad
\vect u=P_u\hat{\vect u},
\label{eq:scaling}
\end{equation}
where $P_r,P_v,P_u$ are positive definite. The implementation uses scalar blocks and $P_u=s_uI_3$, which preserves the spherical thrust set in normalized coordinates. Scaling is applied before condensation because it affects both the terminal row balance and the spectral stepsize \cite{ross2018scaling}.

\section{Exact Dynamics Condensation}

\subsection{Discrete State-Space Form}

Let $\vect x_k=[\vect r_k^\mathsf T,\vect v_k^\mathsf T]^\mathsf T$. Eliminating $\vect v_{k+1}$ from Eq.~\eqref{eq:trapr} yields
\begin{equation}
\vect x_{k+1}=A\vect x_k+B_0\vect u_k+B_1\vect u_{k+1}+\vect d,
\label{eq:ss}
\end{equation}
with
\begin{equation}
A=\begin{bmatrix}I_3&\Delta tI_3\\0&I_3\end{bmatrix},\qquad
B_0=B_1=\begin{bmatrix}\frac{\Delta t^2}{4}I_3\\[1mm]\frac{\Delta t}{2}I_3\end{bmatrix},
\label{eq:AB}
\end{equation}
and
\begin{equation}
\vect d=\begin{bmatrix}\frac{\Delta t^2}{2}\vect g\\[1mm]\Delta t\vect g\end{bmatrix}.
\label{eq:d}
\end{equation}
For $S_x=\operatorname{blkdiag}(P_r,P_v)$ and $S_u=P_u$, define
\begin{equation}
\hat A=S_x^{-1}AS_x,\qquad
\hat B_i=S_x^{-1}B_iS_u,\qquad
\hat{\vect d}=S_x^{-1}\vect d.
\label{eq:normmats}
\end{equation}
This change of variables is invertible and does not alter the trapezoidal discretization.

\subsection{Terminal-Sensitivity Recursion}

Stack the normalized controls as
\begin{equation}
\hat{\vect U}=
[\hat{\vect u}_1^\mathsf T,\ldots,\hat{\vect u}_N^\mathsf T]^\mathsf T
\in\R^{3N},
\end{equation}
and let $E_k\in\R^{3\times3N}$ select the $k$th control block. Define a zero-control trajectory $\hat{\vect x}^{(0)}_k$ and sensitivity matrix $\Gamma_k$ by
\begin{align}
\hat{\vect x}^{(0)}_1&=S_x^{-1}\vect x_0,\qquad \Gamma_1=0,
\label{eq:gamma0}\\
\hat{\vect x}^{(0)}_{k+1}&=\hat A\hat{\vect x}^{(0)}_k+\hat{\vect d},
\label{eq:base_rec}\\
\Gamma_{k+1}&=\hat A\Gamma_k+\hat B_0E_k+\hat B_1E_{k+1}.
\label{eq:gamma_rec}
\end{align}
Induction on Eq.~\eqref{eq:ss} gives
\begin{equation}
\hat{\vect x}_k=\Gamma_k\hat{\vect U}+\hat{\vect x}^{(0)}_k .
\label{eq:affine_state}
\end{equation}
Consequently, the complete dynamics chain and terminal boundary reduce to one six-row equality:
\begin{equation}
G\hat{\vect U}=\hat{\vect x}_f-\hat{\vect x}^{(0)}_N,\qquad
G=\Gamma_N\in\R^{6\times3N}.
\label{eq:terminal}
\end{equation}
Every intermediate state can be recovered by one forward recursion after optimization; none is required in the online decision vector.

\subsection{Row Normalization and Exactness}

Terminal position and velocity rows have different units and magnitudes. Let
\begin{equation}
D=\diag\left(\frac{1}{\norm{G_{1,:}}_2},\ldots,
\frac{1}{\norm{G_{6,:}}_2}\right)
\label{eq:row_scale}
\end{equation}
and define
\begin{equation}
H_c=DG,\qquad
\vect h_c=D\left(\hat{\vect x}_f-\hat{\vect x}^{(0)}_N\right).
\label{eq:condensed}
\end{equation}
Each row of $H_c$ has unit Euclidean norm.

\textbf{Proposition 1.} If every row of $G$ is nonzero, then $H_c\hat{\vect U}=\vect h_c$ is equivalent to the terminal condition $\hat{\vect x}_N=\hat{\vect x}_f$. Condensation changes neither the feasible control set nor the value of any control-only objective.

\textit{Proof.} Equation~\eqref{eq:affine_state} makes the terminal condition equivalent to Eq.~\eqref{eq:terminal}. Under the stated assumption, $D$ is invertible; premultiplying a linear equality by $D$ therefore leaves its solution set unchanged. The objective depends only on $\hat{\vect U}$, so its value is also unchanged. \hfill$\square$

\begin{table}[htbp]
\caption{\label{tab:dimensions}Dimension reduction for $N=30$}
\centering
\begin{tabular}{lrr}
\toprule
Quantity & Full epigraph & Condensed composite \\
\midrule
Position and velocity variables & 180 & 0 \\
Control variables & 90 & 90 \\
Fuel epigraph variables & 30 & 0 \\
Dynamics or terminal multipliers & 174 & 6 \\
Total primal dimension & 300 & 90 \\
\bottomrule
\end{tabular}
\end{table}

The exactness result depends on the absence of state-path constraints. A glide-slope, obstacle, state-triggered, or other nodewise constraint would require retained states, partial condensation, or a denser control-space inequality. Thus, full condensation is most attractive when the terminal dimension is small and the path is constrained primarily through control.

\section{Strong Convex Regularization and Closed-Form Proximal Map}

\subsection{Strong-Convexity Construction}

Multiplying Eq.~\eqref{eq:weighted} by a positive common constant and using $P_u=s_uI_3$ gives the normalized condensed problem
\begin{equation}
\min_{\hat{\vect U}}\quad
\frac12\hat{\vect U}^\mathsf TP\hat{\vect U}
+\sum_{k=1}^{N}q_k\norm{\hat{\vect u}_k}_2
+\ind_{\mathcal B}(\hat{\vect U}),
\label{eq:scaled_obj}
\end{equation}
where
\begin{equation}
P_k=\lambda w_kI_3,\qquad
q_k=\frac{1-\lambda}{2}\bar u\,w_k,\qquad
\bar u=\frac{u_{\max}}{s_u},
\label{eq:Pq}
\end{equation}
and $\mathcal B=\{\hat{\vect U}:\norm{\hat{\vect u}_k}_2\leq\bar u,\ \forall k\}$. Because the trapezoidal endpoint weights are one half,
\begin{equation}
\norm P_2=\lambda,\qquad \lambda_{\min}(P)=\frac{\lambda}{2}.
\label{eq:Peigs}
\end{equation}

\textbf{Proposition 2.} For any $0<\lambda\leq1$, the objective in Eq.~\eqref{eq:scaled_obj} is $\mu$-strongly convex with $\mu=\lambda/2$. If the condensed feasible set is nonempty, the optimal control is unique.

\textit{Proof.} Equation~\eqref{eq:Peigs} implies $P\succeq(\lambda/2)I$. The Euclidean norms and the indicator of $\mathcal B$ are closed convex functions. Adding them to the quadratic term preserves its strong-convexity modulus. Restriction to a nonempty closed convex affine slice preserves uniqueness \cite{boyd2004convex}. \hfill$\square$

The proposition explains why the small energy term has algorithmic value beyond smoothing a plotted thrust profile. It replaces a potentially flat pure-fuel optimal face with one unique minimizer and an explicit curvature floor. It does, however, change the primary objective; the resulting solution is fuel dominant rather than mathematically identical to the $\lambda=0$ optimum.

\subsection{Radial Proximal Operator}

For one node, define
\begin{equation}
g_k(\vect u)=q_k\norm{\vect u}_2+
\ind_{\{\norm{\vect u}_2\leq\bar u\}}(\vect u).
\end{equation}
The primal update requires
\begin{equation}
\prox_{\alpha g_k}(\vect z_k)=
\arg\min_{\norm{\vect u}_2\leq\bar u}
\frac12\norm{\vect u-\vect z_k}_2^2+\alpha q_k\norm{\vect u}_2.
\label{eq:prox_problem}
\end{equation}

\textbf{Proposition 3.} Let $r_k=\norm{\vect z_k}_2$. The unique solution of Eq.~\eqref{eq:prox_problem} is
\begin{equation}
\prox_{\alpha g_k}(\vect z_k)=
\begin{cases}
\dfrac{\hat r_k}{r_k}\vect z_k,&r_k>0,\\[2mm]
\zero,&r_k=0,
\end{cases}
\label{eq:prox_closed}
\end{equation}
where
\begin{equation}
\hat r_k=\min\{\max(r_k-\alpha q_k,0),\bar u\}.
\label{eq:radius}
\end{equation}

\textit{Proof.} If $\vect z_k=0$, the origin is the unique minimizer. Otherwise, for fixed $\norm{\vect u}_2=s$, the Cauchy--Schwarz inequality shows that the quadratic distance is minimized when $\vect u$ is aligned with $\vect z_k$. Substituting $\vect u=s\vect z_k/r_k$ reduces Eq.~\eqref{eq:prox_problem} to
\begin{equation}
\min_{0\leq s\leq\bar u}\quad \frac12(s-r_k)^2+\alpha q_ks .
\end{equation}
The unconstrained minimizer $s=r_k-\alpha q_k$ projected onto $[0,\bar u]$ yields Eq.~\eqref{eq:radius}, and substitution recovers Eq.~\eqref{eq:prox_closed}. \hfill$\square$

The operator performs group soft thresholding and then clips the surviving magnitude. It detects zero-thrust nodes when $r_k\leq\alpha q_k$ and active maximum-thrust nodes when $r_k-\alpha q_k\geq\bar u$. No cone epigraph, scalar auxiliary variable, or nested solver is needed. This simplification follows from the shared radial symmetry of the norm and the concentric ball \cite{parikh2014proximal}; it does not directly extend to a nonzero thrust lower bound or an anisotropic actuator set.

\section{Condensed Extrapolated Primal--Dual Method}

\subsection{Saddle-Point Iteration}

With $g(\hat{\vect U})=\sum_kg_k(\hat{\vect u}_k)$, the online problem is
\begin{equation}
\min_{\hat{\vect U}}\quad
\frac12\hat{\vect U}^\mathsf TP\hat{\vect U}+g(\hat{\vect U})
\quad\text{s.t.}\quad H_c\hat{\vect U}=\vect h_c .
\label{eq:compact}
\end{equation}
Its Lagrangian is
\begin{equation}
\mathcal L(\vect U,\vect\eta)=
\frac12\vect U^\mathsf TP\vect U+g(\vect U)
+\vect\eta^\mathsf T(H_c\vect U-\vect h_c).
\end{equation}
Let $\vect\Xi^j$ and $\vect\eta^j$ denote the extrapolated primal and dual variables. The implemented updates are
\begin{align}
\vect U^{j+1}&=\prox_{\alpha g}\!\left[
\vect\Xi^j-\alpha(P\vect\Xi^j+H_c^\mathsf T\vect\eta^j)
\right],
\label{eq:primal}\\
\vect w^{j+1}&=\vect\eta^j+\beta\left[
H_c(2\vect U^{j+1}-\vect\Xi^j)-\vect h_c\right],
\label{eq:dual}\\
\vect\Xi^{j+1}&=(1-\rho)\vect\Xi^j+\rho\vect U^{j+1},
\label{eq:xi}\\
\vect\eta^{j+1}&=(1-\rho)\vect\eta^j+\rho\vect w^{j+1}.
\label{eq:eta}
\end{align}
The six-dimensional product $H_c^\mathsf T\vect\eta$ broadcasts terminal position and velocity error to every control node in one operation. This replaces the propagation of $6(N-1)$ dynamics multipliers along the full state chain.

\subsection{Spectral Stepsizes}

Choose a primal--dual ratio $\omega=\beta/\alpha>0$ and a safety factor $0<c<1$. The steps satisfy
\begin{equation}
\alpha\left(\norm P_2+\beta\norm{H_c}_2^2\right)=c.
\label{eq:step_identity}
\end{equation}
Substituting $\beta=\omega\alpha$ gives the positive root
\begin{align}
\alpha&=\frac{2c}{\norm P_2+
\sqrt{\norm P_2^2+4c\omega\norm{H_c}_2^2}},
\label{eq:alpha}\\
\beta&=\omega\alpha .
\label{eq:beta}
\end{align}
The condensed norm is computed offline from a $6\times6$ Gram matrix because
\begin{equation}
\norm{H_c}_2^2=\lambda_{\max}(H_cH_c^\mathsf T).
\end{equation}
For the $\lambda=0.02$ case, $\norm{H_c}_2=1.36747$, $\omega=0.15$, $\rho=1.995$, and $c=0.99$, resulting in
\begin{equation}
\alpha=1.843372658,\qquad \beta=0.276505899.
\label{eq:num_steps}
\end{equation}
The ratio and extrapolation were selected by an offline candidate search for the fixed test case. Row normalization must precede this search because it changes both $\norm{H_c}_2$ and the dual scaling.

\begin{figure}[htbp]
\centering
\includegraphics[width=0.96\textwidth]{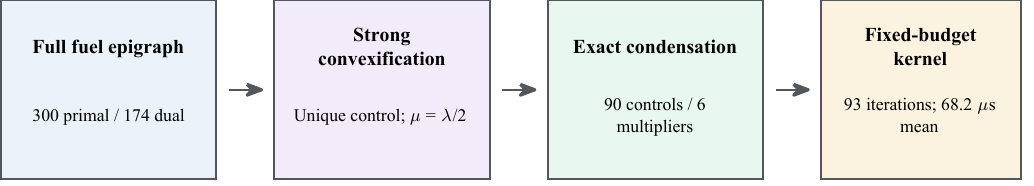}
\caption{\label{fig:pipeline}Structure-exploitation path from a full fuel epigraph to a fixed-budget microsecond solver.}
\end{figure}

\begin{figure}[htbp]
\centering
\fbox{\begin{minipage}{0.92\textwidth}
\small\singlespacing
\textbf{Algorithm 1: Fixed-budget condensed primal--dual iteration}

\textit{Offline:} Construct $H_c$, $P$, and $q$; compute $\norm{H_c}_2$, $\alpha$, and $\beta$; select iteration budget $J_{\max}$.

\textit{Online input:} Current terminal right-hand side $\vect h_c$ and optional warm start $\vect U^0$.

1. Set $\vect\Xi^0=\vect U^0$ and $\vect\eta^0=\zero$.

2. For $j=0,\ldots,J_{\max}-1$:

\quad a. $\vect Z=\vect\Xi^j-\alpha(P\vect\Xi^j+H_c^\mathsf T\vect\eta^j)$.

\quad b. Apply Eq.~\eqref{eq:prox_closed} independently to all $N$ control nodes.

\quad c. Update the six-dimensional dual prediction using Eq.~\eqref{eq:dual}.

\quad d. Extrapolate using Eqs.~\eqref{eq:xi} and \eqref{eq:eta}.

3. Recover the physical control and perform one forward state recursion.
\end{minipage}}
\caption{\label{alg:solver}Online algorithm after all fixed matrices and scalar parameters have been precomputed.}
\end{figure}

\subsection{Residuals and Computational Cost}

The reference-solution metric used in the offline study is
\begin{equation}
e_{\mathrm{ref}}^j=\max\left\{
\frac{\norm{X^j-X^\star}_F}{\max(1,\norm{X^\star}_F)},
\frac{\norm{U^j-U^\star}_F}{\max(1,\norm{U^\star}_F)}
\right\}.
\label{eq:ref_res}
\end{equation}
This metric is unavailable onboard. Implementable alternatives are the equality residual
\begin{equation}
r_{\mathrm{eq}}^j=\norm{H_c\vect U^j-\vect h_c}_\infty
\end{equation}
and the proximal fixed-point residual
\begin{equation}
r_{\mathrm{fp}}^j=
\frac{\norm{\vect U^j-\prox_{\alpha g}[
\vect U^j-\alpha(P\vect U^j+H_c^\mathsf T\vect\eta^j)]}_2}
{\max(1,\norm{\vect U^j}_2)}.
\label{eq:fp_res}
\end{equation}
For the fixed case, $e_{\mathrm{ref}}^j$ first reaches $10^{-3}$ at iteration 93; therefore the timing implementation uses $J_{\max}=93$. A fixed budget removes early-exit branches but is valid only over the offline-tested operating envelope.

Each iteration evaluates one $H_c^\mathsf T\vect\eta$, one $H_c\bar{\vect U}$, one diagonal gradient, and $N$ independent three-dimensional proximal maps. Since $H_c\in\R^{6\times3N}$, the two dense products each require $18N$ scalar multiplications. The per-iteration work and storage are $O(N)$. The full-state sparse formulation is also $O(N)$ when its chain is exploited, so condensation improves dimensions, constants, and data movement rather than the asymptotic order.

\section{Fixed-Array Implementation}

\subsection{Offline--Online Split}

The offline stage constructs $H_c$, computes its spectral norm, generates $P$ and $q$, and verifies the reference solution. When final time, grid, scaling, and dynamics remain fixed, $H_c$ is reusable; a changed navigation state or landing target updates only $\vect h_c$ through Eq.~\eqref{eq:condensed}. In repeated guidance, a shifted previous control can also warm start $\vect U^0$.

Controls are stored node first:
\begin{equation}
[u_{x,1},u_{y,1},u_{z,1},u_{x,2},\ldots].
\end{equation}
This layout makes each proximal triplet contiguous. The six rows of $H_c$ are stored contiguously for sequential access during terminal-residual evaluation. The core C17 routine does not call dynamic memory, a file system, a Basic Linear Algebra Subprogram library, or an optimization library. All array extents are compile-time constants.

\subsection{Memory, Control Flow, and Numerical Protection}

The local solve workspace is approximately \SI{4.4}{\kilo\byte}, including primal, extrapolated, temporary, prediction, dual, and final-recovery arrays. The 540 double-precision constants in $H_c$ require \SI{4.32}{\kilo\byte}. Reference controls and reference states are used only for offline validation and can be removed from a deployment build.

The real-time loop always executes 93 iterations and contains no matrix factorization, line search, or variable-length inner problem. When $r_k=0$, the proximal implementation sets its scale explicitly to zero to prevent division by zero. The reported build uses IEEE 754 double precision, Microsoft Visual C++ 2022, release optimization, whole-program optimization, and precise floating-point semantics. Flight software would additionally require finite-input checks, floating-point exception handling, output guards, and a defined failure-management path.

\section{Numerical Results}

\subsection{Scenario and Reference Solutions}

The initial state and gravity are
\begin{align}
\vect r_0&=[500,2500,-2000]^\mathsf T\ \si{\meter},\\
\vect v_0&=[0,-75,-100]^\mathsf T\ \si{\meter\per\second},\\
\vect g&=[0,-3.71,0]^\mathsf T\ \si{\meter\per\second\squared}.
\end{align}
The target position and velocity are zero, $t_f=\SI{60}{\second}$, and $N=30$. The modeled mass is \SI{1905}{\kilogram}; maximum thrust is \SI{19200}{\newton}, giving $u_{\max}=\SI{10.0787}{\meter\per\second\squared}$. The scale matrices are $P_r=0.943263I_3$, $P_v=0.843542I_3$, and $P_u=1.256780I_3$.

The pure-fuel reference is obtained from the full conic epigraph using the MATLAB \texttt{coneprog} routine. The weighted reference uses a high-accuracy condensed primal--dual warm start followed by sequential quadratic programming with $10^{-10}$ constraint and optimality tolerances. Reference generation is excluded from all reported online timings. Every method uses the same grid, boundary conditions, thrust limit, and scaling.

\subsection{Trajectory Accuracy}

\begin{figure}[htbp]
\centering
\includegraphics[width=0.98\textwidth]{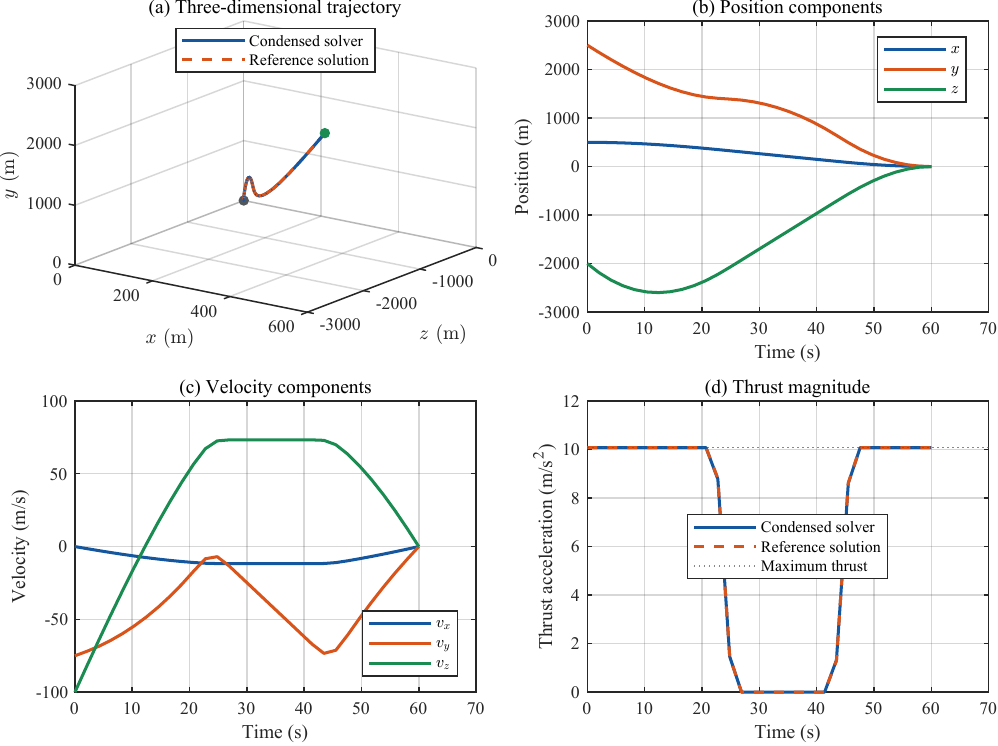}
\caption{\label{fig:trajectory}Landing trajectory and thrust history for the fuel-dominant setting $\lambda=0.02$.}
\end{figure}

Figure~\ref{fig:trajectory} compares the 93-iteration solution with the weighted reference. The trajectories and thrust magnitudes are visually coincident, and the thrust bound is respected at every node. The terminal state after 93 iterations is approximately
\begin{align}
\vect r_N&=[0.0051,-0.6353,-0.2649]^\mathsf T\ \si{\meter},\\
\vect v_N&=[-0.00039,0.01586,0.00929]^\mathsf T\
\si{\meter\per\second}.
\end{align}
These errors correspond to the common relative comparison threshold in Eq.~\eqref{eq:ref_res}; they are not asserted to be mission-level terminal tolerances. Tighter physical tolerances require a larger fixed budget or direct residual-based stopping.

\subsection{Convergence and Structural Reduction}

\begin{figure}[htbp]
\centering
\includegraphics[width=0.98\textwidth]{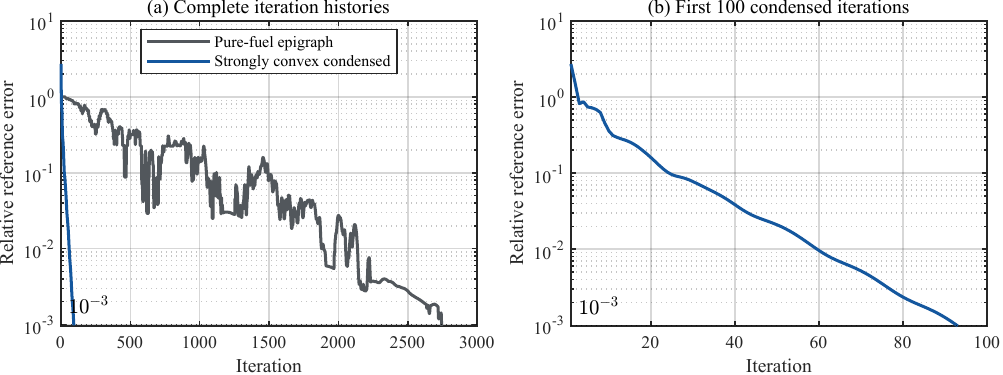}
\caption{\label{fig:convergence}Convergence of the full pure-fuel epigraph and the strongly convex condensed formulation.}
\end{figure}

The full pure-fuel epigraph and strongly convex condensed formulations reach $e_{\mathrm{ref}}\leq10^{-3}$ in 2744 and 93 iterations, respectively, a factor of 29.5 in iteration count. Three mechanisms act together: positive quadratic curvature suppresses drift along a nearly flat pure-fuel optimum set; $H_c^\mathsf T\vect\eta$ communicates terminal error directly to all nodes; and the radial proximal map removes the scalar fuel epigraph. The available experiment measures their combined effect and should not be interpreted as a single-factor ablation.

\begin{table}[htbp]
\caption{\label{tab:efficiency}Numerical and computational comparison}
\centering
\small
\begin{tabular}{p{0.24\textwidth}p{0.19\textwidth}rrrr}
\toprule
Method & Objective and structure & Primal/dual & Iter. & $J_F$, \si{\meter\per\second} & Mean time \\
\midrule
Full epigraph (MATLAB) & Pure fuel; full state & 300/174 & 2744 & 396.1060 & \SI{310.3}{\milli\second} \\
Condensed (MATLAB) & $\lambda=0.02$; condensed & 90/6 & 93 & 396.2356 & \SI{34.9}{\milli\second} \\
Validation mode (C) & Condensed; validation & 90/6 & 93 & 396.2356 & \SI{108.6}{\micro\second} \\
Fixed-budget mode (C) & Condensed; timing & 90/6 & 93 & 396.2356 & \SI{68.2}{\micro\second} \\
\bottomrule
\end{tabular}
\end{table}

\subsection{Curvature--Performance Tradeoff}

\begin{figure}[htbp]
\centering
\includegraphics[width=0.98\textwidth]{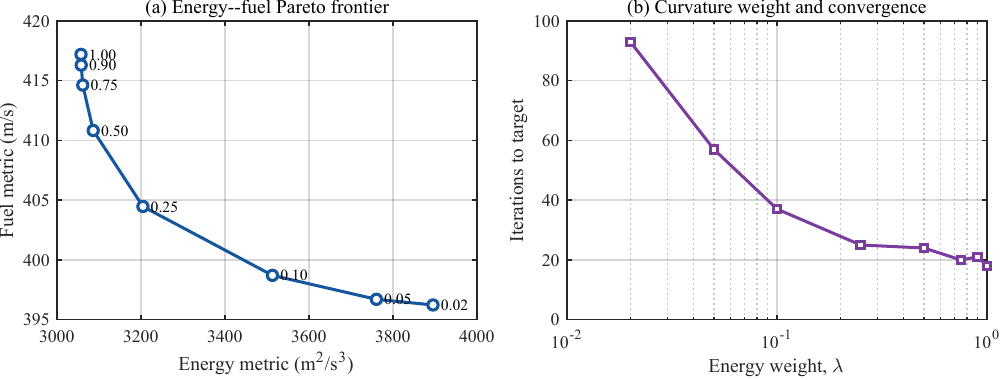}
\caption{\label{fig:weight}Energy weight governs the energy--fuel tradeoff and the curvature available to the first-order iteration.}
\end{figure}

Reference solutions were generated for $\lambda\in\{0.02,0.05,0.10,0.25,0.50,0.75,0.90,1\}$, with stepsize parameters selected from a common offline candidate set. As $\lambda$ increases from 0.02 to 1, the strong-convexity modulus increases from 0.01 to 0.5 and the best observed iteration count generally decreases from 93 to 18. Energy consumption decreases while the fuel metric increases, forming a clear Pareto frontier.

\begin{table}[htbp]
\caption{\label{tab:weights}Effect of the energy weight}
\centering
\small
\begin{tabular}{rrrrrrr}
\toprule
$\lambda$ & $\mu$ & Iter. & $\omega$ & $\rho$ & $J_E$, \si{\meter\squared\per\second\cubed} & $J_F$, \si{\meter\per\second} \\
\midrule
0.02 & 0.010 & 93 & 0.15 & 1.995 & 3894.89 & 396.23 \\
0.05 & 0.025 & 57 & 0.30 & 1.995 & 3760.08 & 396.70 \\
0.10 & 0.050 & 37 & 1.00 & 1.980 & 3512.58 & 398.71 \\
0.25 & 0.125 & 25 & 1.00 & 1.980 & 3204.21 & 404.47 \\
0.50 & 0.250 & 24 & 1.50 & 1.950 & 3085.85 & 410.83 \\
0.75 & 0.375 & 20 & 3.00 & 1.950 & 3060.94 & 414.64 \\
0.90 & 0.450 & 21 & 3.00 & 1.900 & 3057.24 & 416.30 \\
1.00 & 0.500 & 18 & 5.00 & 1.950 & 3056.75 & 417.20 \\
\bottomrule
\end{tabular}
\end{table}

\begin{figure}[htbp]
\centering
\includegraphics[width=0.76\textwidth]{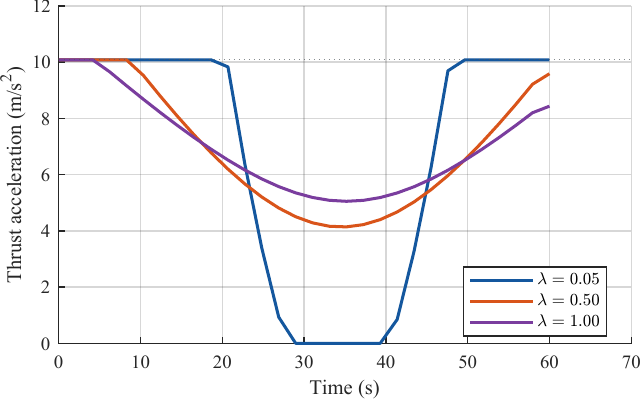}
\caption{\label{fig:controls}Reference thrust magnitudes for representative strong-convexity weights.}
\end{figure}

The 93-iteration $\lambda=0.02$ solution has $J_F=\SI{396.2356}{\meter\per\second}$, an increase of \SI{0.1296}{\meter\per\second}, or 0.033\%, relative to the pure-fuel value. The corresponding weighted reference has $J_F=\SI{396.2316}{\meter\per\second}$, showing that most of the reported difference is the deliberate objective tradeoff rather than incomplete convergence. At small $\lambda$, group shrinkage also creates 20--26.7\% zero-thrust nodes; larger quadratic weights distribute control over more nodes, as shown in Fig.~\ref{fig:controls}.

\subsection{Cross-Language Agreement and Timing}

The fixed-array validation mode and the high-level implementation use the same constants, initialization, and 93 updates. The maximum absolute discrepancies in position, velocity, and control are $4.41\times10^{-12}$~m, $1.14\times10^{-13}$~m/s, and $3.02\times10^{-14}$~m/s$^2$, respectively. Energy and fuel metrics agree to the displayed precision. This test covers array layout, trapezoidal state recovery, all proximal branches, the primal--dual update, and conversion back to physical units.

\begin{figure}[htbp]
\centering
\includegraphics[width=0.74\textwidth]{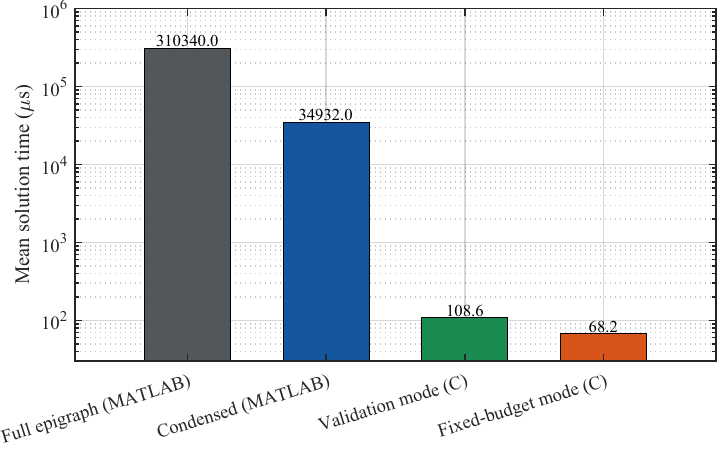}
\caption{\label{fig:runtime}Mean end-to-end solve time; the vertical axis is logarithmic.}
\end{figure}

The height of each bar in Fig.~\ref{fig:runtime} is a measured mean
wall-clock time for one end-to-end solve, expressed in microseconds. From left
to right, the bars have the following meanings. The \emph{full epigraph
(MATLAB)} bar, \SI{310.340}{\milli\second}, is the 2744-iteration pure-fuel
baseline with 300 primal variables, 174 dynamics multipliers, and the full
state trajectory retained. The \emph{condensed (MATLAB)} bar,
\SI{34.932}{\milli\second}, is the strongly convex $\lambda=0.02$ problem
with 90 controls, six terminal multipliers, 93 updates, and high-level
history and reference-error bookkeeping. The \emph{validation mode (C)}
bar, \SI{108.6}{\micro\second}, performs the same 93 numerical updates in
fixed arrays but reconstructs the trajectory and compares it with the
stored reference solution at every iteration; it is intended to establish
the first iteration that reaches the $10^{-3}$ threshold. The
\emph{fixed-budget mode (C)} bar, \SI{68.2}{\micro\second}, removes those
per-iteration state reconstructions and reference comparisons, always
executes 93 updates, and performs the terminal residual, one trajectory
reconstruction, reference check, and energy and fuel calculations only
after the iteration loop. Thus, the third-to-fourth difference isolates
validation bookkeeping and repeated state recovery; larger gaps to the first
two bars also include formulation, iteration count, runtime, memory layout,
history collection, and compiler effects.

Timing was performed on an Intel Core i7-10875H processor with Microsoft
Visual C++ 2022, release x64 optimization, and one calling thread. After 200
warmup solves, 10,000 complete fixed-budget solves were timed with
\texttt{QueryPerformanceCounter}. The mean was
\SI{68.2}{\micro\second}; observed 95th-percentile values ranged from 90.1
to \SI{103.6}{\micro\second}, and the 99th percentile was
\SI{128.1}{\micro\second}. Occasional Windows scheduling outliers
approached \SI{1.20}{\milli\second}. The mean is approximately 512 times
smaller than the \SI{34.9}{\milli\second} high-level condensed
implementation on the same computer. Because the two implementations have
different bookkeeping and execution environments, this ratio is an
end-to-end implementation comparison, not a language-only or
formulation-only speedup.

\section{Discussion}

\subsection{Why the Computation Reaches the Microsecond Regime}

The timing result follows from a sequence of formulation decisions, not from a single low-level optimization. Strong convex regularization supplies a unique target and positive curvature; exact condensation shrinks feedback from 174 dynamics multipliers to six terminal multipliers; the joint proximal map eliminates 30 epigraph variables and all cone-projection subiterations; spectral precomputation removes online norm estimation; and a fixed budget makes the execution path deterministic at the algorithmic level. The online kernel therefore consists almost entirely of streaming multiply--accumulate operations, three-dimensional norms, scalar comparisons, and vector scaling.

The numerical evidence does not isolate the acceleration attributable to each component. A controlled ablation would hold the objective, target tolerance, parameter-search budget, implementation language, and timing boundary fixed while independently toggling state condensation, epigraph elimination, strong convexification, and warm starting. The present results support the combined architecture and its mechanism, not an unsupported percentage allocation among components.

\subsection{Fuel Optimality and Convexification Scope}

For $\lambda>0$, the solved problem is the strongly convex energy--fuel objective in Eq.~\eqref{eq:weighted}, not the strict $\lambda=0$ pure-fuel problem. The 0.033\% increase is small in the present case, but it remains a measurable objective perturbation. Accordingly, the method is described as fuel dominant or near-fuel-optimal when discussing $\lambda=0.02$.

The current actuator set has zero minimum thrust and is convex before optimization. Standard lossless convexification can produce a convex subproblem for a nonzero lower-thrust bound under additional assumptions, but that formulation introduces a different feasible-set geometry and generally destroys the concentric-ball proximal formula used here. Extending the microsecond architecture to nonzero minimum thrust, pointing constraints, mass depletion, or six-degree-of-freedom motion requires a new proximal decomposition or use of this condensed solver inside a successive-convexification loop.

\subsection{From Soft Real Time to Flight Qualification}

Mean and percentile timing on a desktop Windows system demonstrate a small computational workload but do not constitute a worst-case execution-time certificate. Flight deployment requires target-processor measurements, fixed clock and cache policies, interrupt control, compiler qualification, stack analysis, worst-case execution-time analysis, hardware-in-the-loop testing, and Monte Carlo coverage over navigation errors, target changes, and thrust uncertainty.

The fixed $H_c$ matrix can be reused only while the grid, final time, scaling, and linear dynamics are unchanged. Receding-horizon guidance with a changed initial state can update $\vect h_c$ cheaply, but free-final-time or nonlinear models require matrix reconstruction and new stepsize bounds. Path constraints also reduce the benefit of full condensation. These limitations define the intended result precisely: a microsecond-level convex optimization kernel for terminally constrained translational landing, rather than a complete flight-ready guidance system.

\section{Conclusions}

This paper presents a direct path from a convex powered-descent model to fixed-size C code. A low-weight energy term gives the fuel-dominant objective strong-convexity curvature. Terminal-sensitivity recursion reduces the full-state epigraph with 300 primal variables and 174 equality multipliers to 90 controls and six terminal equalities. The shared radial geometry of the fuel norm and thrust ball gives an exact closed-form proximal operator. Each condensed update therefore uses two $6\times90$ matrix--vector products, a diagonal gradient, and 30 independent three-dimensional proximal evaluations.

At $\lambda=0.02$, xPIPG reaches a relative reference-solution error of $10^{-3}$ in 93 updates, compared with 2744 updates for the pure-fuel full-state epigraph baseline. The fuel metric increases by 0.033\%. The C17 and MATLAB controls agree to double-precision rounding, and the mean end-to-end solve time on the tested desktop processor is \SI{68.2}{\micro\second}.

The result demonstrates the combined value of strong convex regularization, exact condensation, an exact radial proximal map, and fixed-budget implementation. It is not a new lossless-convexification proof for a general nonconvex powered-descent problem, nor is it a worst-case timing certificate for a flight processor. Future work will address warm starts over changing initial conditions, controlled ablations, target-hardware timing, mass depletion, nonzero lower thrust, pointing limits, and six-degree-of-freedom dynamics.

\section*{Funding Sources}

This research received no external funding.

\section*{Declarations}

None.

\bibliography{references}

\end{document}